\documentclass[12pt]{amsart}

\usepackage{amsmath,amsthm,amssymb,amsfonts,bbm,latexsym,xcolor,mathrsfs,hyperref,comment,soul}

\allowdisplaybreaks

\numberwithin{equation}{section}

\renewcommand*{\theequation}{%
  \ifnum\value{subsection}=0 %
    \thesection
  \else
    \thesubsection
  \fi
  .\arabic{equation}%
}

\title[Multiple Dirichlet series predictions for moments of $L$-functions]{Multiple Dirichlet series predictions for moments of $L$-functions: unitary, symplectic and orthogonal examples}
\author{Siegfred Baluyot}
\address{Mathematics Department\\ East Carolina University \\ Greenville, NC 27858}
\email{baluyots24@ecu.edu}
\author{Martin \v{C}ech}
\address{Charles University, Faculty of Mathematics and Physics, Department of Mathematical Analysis and Department of Algebra, Sokolovska´ 83, 18600 Praha 8, Czech Republic}
\email{\href{mailto:martin.cech@matfyz.cuni.cz }{martin.cech@matfyz.cuni.cz}}


\newcommand{\lz}{\left(}
\newcommand{\pz}{\right)}

\newcommand{\C}{\mathbb{C}}
\newcommand{\Z}{\mathbb{Z}}
\newcommand{\sumstar}{\sideset{}{^{*}}\sum}

\newcommand{\R}{\mathbb{R}}

\newcommand{\lab}{\left|}
\newcommand{\rab}{\right|}

\newcommand{\re}{\mathrm{Re}}

\newcommand{\bfrac}[2]{\lz\frac{#1}{#2}\pz}

\newcommand{\res}[1]{\underset{#1}{\mathrm{Res\ }}}
\renewcommand{\mod}[1]{\text{ (mod $#1$)}}

\renewcommand{\phi}{\varphi}

\newcommand{\sumflat}{\sideset{}{^{\flat}}\sum}

\newtheorem{theorem}{Theorem}[section]
\newtheorem{conjecture}[theorem]{Conjecture}

\newtheorem{prop}[theorem]{Proposition}

\newtheoremstyle{remarks}%
{3pt}
{3pt}
{}
{}
{\bfseries}
{.}
{5pt}
{}%

\theoremstyle{remarks}
\newtheorem{remark}[theorem]{Remark}

\begin{document}

\begin{abstract}
We devise heuristics using multiple Dirichlet series to predict asymptotic formulas for shifted moments of (1) the family of Dirichlet $L$-functions of all even primitive characters of conductor $\leq Q$, with $Q$ a parameter tending to infinity, (2) the family of quadratic Dirichlet $L$-functions, (3) the family of quadratic twists of an $L$-function associated to a fixed Hecke eigencuspform for the full modular group, and (4) the family of quadratic twists of an $L$-function of a fixed arbitrary elliptic curve over $\mathbb{Q}$ that has a non-square conductor. For each of these families, the resulting predictions agree with the predictions of the recipe developed by Conrey, Farmer, Keating, Rubinstein, and Snaith, except for (4), 
 where the recipe requires a slight modification due to a correlation between the Dirichlet coefficients and the root number of the corresponding $L$-functions. We find a one-to-one correspondence between the residues from the multiple Dirichlet series analysis and the terms from the recipe prediction.

\end{abstract}

\maketitle


\section{Introduction}

Understanding the behavior of moments of $L$-functions is an important problem in analytic number theory. In the case of the Riemann zeta-function, asymptotic formulas are known only for the second and fourth moments. The situation is similar for other families of $L$-functions, where asymptotic formulas are also known only for a few low moments. In the late 90's, Keating and Snaith~\cite{KS00a,KS00b,KS03} modeled $L$-functions using characteristic polynomials of large random matrices. Their heuristic is the first in history that results in viable predictions for the exact values of the leading terms of \textit{all} moments of various families of $L$-functions. Since then, other heuristic methods for predicting asymptotic formulas have been developed.

One of these heuristics is the \textit{recipe} devised by Conrey, Farmer, Keating, Rubinstein, and Snaith (CFKRS)~\cite{cfkrs}. One feature of the recipe is that it is readily applicable to a wide variety of families of $L$-functions. In fact, CFKRS formulated it for a general family of $L$-functions. CFKRS give examples of predictions for each of the three symmetry types of families of $L$-functions: unitary, symplectic, and orthogonal in the context of the Katz-Sarnak philosophy~\cite{katzsarnak}. The predictions of the CFKRS recipe are consistent with the predictions of Keating and Snaith, even though the recipe does not use any random matrix theory. 

We remark that the individual steps in the recipe ``cannot be justified'' (see Note (6) on Page 21 of \cite{cfkrs}) because some terms that are discarded or added back in can be of the same size as the main term. CFKRS conjecture that these large errors cancel and the recipe provides a correct prediction. Moreover, in practice, asymptotic formulas can arise in a different form than the recipe prediction and it can take some effort to show that the formulas are in agreement (see, for instance, \cite{btb} or \cite{CoRo}).

\smallskip

Another heuristic for predicting asymptotic formulas comes from using multiple Dirichlet series. This heuristic was first developed by Diaconu, Goldfeld, and Hoffstein~\cite{dgh}, who studied the multiple Dirichlet series associated with integral moments of the Riemann zeta function and the family of quadratic Dirichlet $L$-functions. They formulated natural conjectures about the poles, residues, and meromorphic continuation of the multiple Dirichlet series. Their conjectures, if true, would imply the Keating--Snaith prediction. The most challenging part of the method that prevents us from evaluating high moments in these families is in showing there is a meromorphic continuation of the multiple Dirichlet series to a wide enough region.

Multiple Dirichlet series have also been useful in proving asymptotic formulas for low moments of quadratic Dirichlet $L$-functions~\cite{bucurdiaconu,cech,diaconu19,dgh,DW21,goho,HR92}, quadratic twists of $L$-functions~\cite{BFH2,FF04}, and $L$-functions in the $t$-aspect~\cite{good86,zhang05,zhang06}.

\smallskip

In this paper, we further develop the multiple Dirichlet series heuristic and use it to predict asymptotic formulas for shifted moments of four families of $L$-functions of three different symmetry types. Our heuristic provides a direct comparison between the recipe and the multiple Dirichlet series prediction. In particular, we will note that there is a one-to-one correspondence between the terms arising in both predictions. We will also state conjectures about the analytic properties of the multiple Dirichlet series that imply the asymptotic formula for the moments precisely in the form predicted by the recipe. We consider the following families:
\begin{enumerate}
    \item the family of Dirichlet $L$-functions $L(s,\chi)$ with $\chi$ varying over all even primitive Dirichlet characters of conductor $\leq Q$, where $Q$ is a parameter tending to $\infty$,
    \item the family of quadratic Dirichlet $L$-functions $L(s,\chi_d)$, where $\chi_d= (d|\cdot)$ is the Kronecker symbol and $d$ varies over all positive fundamental discriminants $\leq X$, with $X$ a parameter tending to $\infty$,
    \item the family of twisted $L$-functions $L(s,f\otimes \chi_d)$, where $f$ is an arbitrary fixed Hecke eigencuspform for the full modular group $\mathrm{SL}(2,\Z)$, where $d$ varies over all positive fundamental discriminants $\leq X$, and
    \item the family of twisted $L$-functions $L(s,E_d)$, where $E$ is an arbitrary fixed elliptic curve over $\mathbb{Q}$ with a non-square conductor, $E_d$ is the quadratic twist of $E$ by $d$, and $d$ varies over all positive fundamental discriminants $\leq X$.
\end{enumerate}
For each of these families, we provide the CFKRS recipe prediction, and then study the associated multiple Dirichlet series. We find its poles and residues in a region of absolute convergence. We also find its functional equations that arise from the functional equations of the individual $L$-functions, and note that they provide new poles of the multiple Dirichlet series. We then state a conjecture about the meromorphic continuation and its properties, apply Perron's formula, shift the resulting integral and capture the contribution of the residues to arrive at a (conjectural) asymptotic formula. For the first three of these families, the resulting predictions agree with the predictions of the CFKRS recipe. On the other hand, in Case (4), we find a disagreement between the two predictions. This discrepancy is caused by the correlation between the root numbers and the Dirichlet series coefficients of the functions $L(s,E_d)$. For completeness, we also describe how the CFKRS recipe reaches its predictions for each of these four families. Furthermore, we suggest a slight modification of the recipe for Case (4) so that the prediction of the modified version agrees with the prediction of the multiple Dirichlet series heuristic. In the Appendix, we verify that the work of Shen~\cite{shen} shows that the multiple Dirichlet series prediction, and thus the prediction of the modified recipe, is correct for the first moment of the family.

In comparing the two heuristics, our main observation is that there is a one-to-one correspondence between the terms in the recipe prediction and the residues in the multiple Dirichlet series prediction. This is true for Cases (1), (2), and (3), and is true for Case (4) after modifying the recipe. More precisely, for each subset $T\subset\{\alpha_1,\dots,\alpha_k\}$ of the set of shift variables in the shifted moment, the term from the recipe prediction that results from using the second part of the approximate functional equation for the $L$-functions $L(\frac{1}{2}+\alpha_j,f)$ with $\alpha_j\in T$ is exactly the same as the multiple Dirichlet series residue arising from using the functional equation for the same $L$-functions. This correspondence is perhaps not too surprising, as the approximate functional equation may be considered as a reformulation of the functional equation. Nevertheless, it provides a clear connection between the two approaches and gives a more transparent interpretation of some of the results and computations of Diaconu, Goldfeld, and Hoffstein~\cite{dgh}, which we also extend to other families.

\smallskip

The multiple Dirichlet series approach to moments of $L$-functions has been the subject of a great deal of research in recent years. Goldfeld and Hoffstein~\cite{goho} used this approach to improve the asymptotic formula found earlier by Jutila~\cite{jutila} for the first moment of quadratic Dirichlet $L$-functions. Hoffstein and Rosen~\cite{HR92} have deduced analogues of the results of Goldfeld and Hoffstein for function fields over finite fields. In the case of this family in the number field setting, Diaconu and Whitehead~\cite{DW21} have proved an asymptotic formula for the (unshifted) third moment with error term of size at most $2/3+\varepsilon$ the power of the leading order term. In doing so, they have discovered the existence of a term of size about $3/4$ the power of the leading order term in the asymptotic formula. Diaconu~\cite{diaconu19} has proved an analogous result in the function field setting, and found the analogous term in the asymptotic formula. Such lower order terms have not been found in prior work that use the more standard approach of using approximate formulas (e.g.,  approximate functional equations) for the $L$-functions~\cite{florea17,sound00,young13}, and have also not been predicted by the CFKRS recipe or by random matrix theory heuristics. More recently, Diaconu and Twiss~\cite{DT23} have used multiple Dirichlet series techniques to predict all the secondary terms in the asymptotic formula for all (unshifted) integer moments of the family of quadratic Dirichlet $L$-functions in the function field setting, up to a square root cancellation error term.

The multiple Dirichlet series approach has also been useful in studying quadratic twists of $L$-functions~\cite{BFH2,cech2,dgh,FF04,FH95,zhang06B}. One reason that this approach is effective in investigating families of quadratic Dirichlet $L$-functions and families of quadratic twists of $L$-functions is that these families are parametrized by integers and thus it is natural to apply Perron's formula (or Mellin inversion) to relate sums over the family to integrals that involve multiple Dirichlet series.  Moreover, these multiple Dirichlet series have (at least heuristically) an extra functional equation, which should be thought of as analogous to Poisson summation (the ``harmonic detector" in the classical setting).  This extra functional equation is not important in this paper, but plays a major role in extending the regions of meromorphic continuation (see, for example, the surveys \cite{bump} and \cite{BFH2}). Our first family studied in Section \ref{sec: unitary MDS} is not a family of quadratic twists and has not been studied before using multiple Dirichlet series. The remaining three families have the form of quadratic twists of a fixed $L$-function.

\smallskip

We now outline the rest of the paper. In Section~\ref{sec: notations}, we describe notations that are used throughout the paper. We discuss the family of Dirichlet $L$-functions of even primitive characters in Section~\ref{sec: unitary}. The discussion of this family is more detailed and includes several remarks that similarly apply to the other families. We discuss the family of quadratic Dirichlet $L$-functions in Section~\ref{sec: symplectic}, the family of quadratic twists of $L$-functions associated to a Hecke eigencuspform for SL$(2,\mathbb{Z})$ in Section~\ref{sec: orthogonal}, and the family of $L$-functions associated to quadratic twists of an elliptic curve in Section~\ref{sec: elliptic}. For each of these families, we present
both the CFKRS heuristic, which we follow directly from \cite{cfkrs}, and the multiple Dirichlet series heuristic we devise (based on earlier works \cite{cech}, \cite{cech2}, \cite{dgh}). In the Appendix, we explain how to modify the argument of Shen~\cite{shen} to deduce an asymptotic formula for the first moment of the family of $L$-functions associated to quadratic twists of an elliptic curve  and show that the (unmodified) recipe leads to a wrong answer for this family, at least for the first moment.

We also remark that similar results and predictions would hold for other families of $L$-functions, and other related problems, for instance the Ratios conjectures (see \cite{cosn} for the Ratios conjecture and its applications; a specific case using multiple Dirichlet series was worked out in \cite{cech}) or the one-level density of low-lying zeros in families of $L$-functions.

\subsection{Notations}\label{sec: notations}
Throughout this paper, $g(x)$ is a smooth, compactly supported test function, and we use $\tilde g(s)$ to denote its Mellin transform, which is defined by
$$
\tilde g(s)=\int_0^\infty g(x)x^{s-1} dx.
$$
Since $g(x)$ is smooth and compactly supported, repeated integration by parts shows that $\tilde g(s)$ has a meromorphic continuation to all of $\mathbb{C}$, with possible poles only at the non-positive integers. The same reasoning implies that in any fixed vertical strip, $\tilde g(s)$ decays faster than any polynomial, i.e.,
$$|\tilde g(\sigma+it)|\ll_{A } (1+|t|)^{-A}$$ as $|t|\rightarrow \infty$ for any $A>0$, uniformly for $\sigma$ in any fixed closed interval.

We let $S=\{s_1,\dots,s_k\}$ denote a multiset of complex numbers. For $J\subseteq\{1,\dots,k\}$, we define
\begin{equation}\label{eqn: sjjdef}
S_J:=\{s_j: j\in J\}, \hspace{15pt} S_J^-:=\{1-s_j: j\in J\}, \hspace{15pt} s_j^J=\begin{cases}
1-s_j &\text{if $j\in J$,}\\
s_j &\text{if $j\notin J$.}
\end{cases}
\end{equation}
We also let $Z=\{z_1,\dots,z_k\}$ and define  $Z_J$, $Z_J^-$, and $z_j^J$ similarly as in \eqref{eqn: sjjdef}.

We define a \textit{tube domain} in $\C^n$ to be a set $T$ of the form $$T=V+i\R^n,$$ where $V\subseteq\R^n$ is open and connected. Tube domains are multidimensional analogues of vertical strips in $\mathbb{C}$. This makes them natural regions of definition of multiple Dirichlet series.

\section*{Acknowledgments}
This work began during the second author's visit at the American Institute of Mathematics, and both authors are grateful for the Institute's hospitality and to Brian Conrey for encouragement and some useful discussions. The first author was supported by NSF DMS-1854398 FRG. The second author was supported by bourse de doctorat en recherche (B2X) of Fonds de recherche
du Qu\'ebec – Nature et technologies (FRQNT), the Research Council of Finland grant no. 333707, and the Charles University Primus programme PRIMUS/24/SCI/010 and PRIMUS/25/SCI/017.

\section{Unitary family: Dirichlet $L$-functions of even primitive characters of conductor $\leq Q$} \label{sec: unitary}
In this section, we consider the family of $L$-functions $L(s,\chi)$ of even primitive Dirichlet characters $\chi$ modulo $q$, averaged over $q$. These $L$-functions satisfy the functional equation 
\begin{equation}\label{eqn: fe unitary}
L(s,\chi)=q^{1/2-s}\varepsilon(\chi)\mathcal{X}(s) L(1-s,\overline\chi),
\end{equation}
where
\begin{equation}\label{eqn: Xdef}
\mathcal{X}(s) = \pi^{s-\frac{1}{2}}\frac{\Gamma\bfrac{1-s}{2}}{\Gamma\bfrac s2},
\end{equation}
and $\varepsilon(\chi)$ is the root number satisfying $|\varepsilon(\chi)|=1$.
We also have the approximate functional equation
\begin{equation}\label{eqn: approximate fe unitary}
L(s,\chi) \approx \sum_n \frac{\chi(n)}{n^s} + \varepsilon(\chi) q^{1/2-s} \mathcal{X}(s) \sum_n \frac{\overline{\chi}(n)}{n^{1-s}}.
\end{equation}

\begin{remark}
    In the recipe, after some steps that include discarding the off-diagonal terms, we eventually extend the $n$-sums in \eqref{eqn: approximate fe unitary} to infinity. For this reason, it is not necessary to indicate the range of $n$ in the $n$-sums. A precise approximate functional equation (see, for instance, \cite[Theorem 5.3]{IK}) usually includes a weight function or an error term.
\end{remark}

We wish to understand the twisted shifted moment 
\begin{align}
\mathcal{U}_{k}(Q;S,Z;M,N)=\sum_{q\geq 1}g\bfrac qQ\sumflat_{\chi\mod q} \overline{\chi}(M) \chi (N) L(s_1,\chi)\dots L(s_k,\chi)L(z_1,\overline\chi)\dots L(z_k,\overline\chi), \label{eqn: unitarymoment}
\end{align}
where $Q$ is a parameter tending to $\infty$, $k,M,N$ are positive integers, and $S=\{s_1,\dots,s_k\}$ and $Z=\{z_1,\dots,z_k\}$ are multisets of complex numbers close to $1/2$, say $s_j-\frac{1}{2}\ll (\log Q)^{-1}$ and $z_j-\frac{1}{2}\ll (\log Q)^{-1}$ for each $j$. The summation symbol $\sum^{\flat}$ here means that we are summing over all even primitive Dirichlet characters modulo $q$. More generally, we may take $S$ and $Z$ to have different cardinalities, but we assume that $|S|=|Z|=k$ for simplicity. Also, our restriction to even primitive characters is only for convenience, and the procedures in this paper also apply to the odd characters with minor modifications.

\subsection{Unitary family: carrying out the CFKRS recipe} 
We now carry out the recipe described in \cite[Section 4.1]{cfkrs}. Without the average over $q\leq Q$, this family is studied in \cite[Section 4.3]{cfkrs}. We replace each $L(s,\chi)$ in \eqref{eqn: unitarymoment} by its approximate functional equation \eqref{eqn: approximate fe unitary}, and multiply out the resulting product. We then throw away all the terms except those that have the same number of $\varepsilon(\chi)$ factors as $\varepsilon(\overline{\chi})$ factors, as otherwise we expect the average of the root numbers over the family to be essentially $0$, and use the fact that $\varepsilon(\chi)\varepsilon(\overline{\chi})=1$. The sum of the remaining terms is
\begin{align*}
\sum_{\substack{J,H \subseteq \{1,\dots,k\} \\ |J|=|H|}} \sum_{q\geq 1}g\bfrac qQ\sumflat_{\chi\mod q} \overline{\chi}(M) \chi (N) \prod_{j\in J}\frac{\mathcal{X}(s_j)}{q^{s_j-\frac{1}{2}}} \prod_{h\in H} \frac{\mathcal{X}(z_h)}{q^{z_h-\frac{1}{2}}} \\
\times \sum_{\substack{n_1,\dots,n_k \\ m_1,\dots,m_k}} \frac{ \chi\left( \prod_{j\not\in J} n_j \prod_{h\in H} m_h \right) \overline{\chi}\left( \prod_{j\in J} n_j \prod_{h\not\in H} m_h \right) }{n_1^{s^J_1} \cdots n_k^{s^J_k} m_1^{z^H_1} \cdots m_k^{z^H_k} }.
\end{align*}
Now we average over the family, bringing the $q$- and $\chi$-sums inside and use the approximation 
$$
\begin{aligned}
\sum_{q\geq 1}g\bfrac{q}{Q}&\sumflat_{\chi\mod q} \chi(Nn)\overline{\chi}(Mm) \\
&\approx
\begin{cases}
    \sum\limits_{q\geq 1}g\bfrac{q}{Q}\sumflat\limits_{\chi\mod q} 1, &\hbox{if $Nn=Mm$ and $(MNmn,q)=1$,}\\
    \\
    0, & \hbox{otherwise}
\end{cases}
\end{aligned}
$$
to arrive at the prediction
\begin{align}
\mathcal{U}_k(Q;S,Z;M,N) \sim \sum_{\substack{J,H \subseteq \{1,\dots,k\} \\ |J|=|H|}} \sum_{q\geq 1}g\bfrac qQ\sumflat_{\chi\mod q} \prod_{j\in J}\frac{\mathcal{X}(s_j)}{q^{s_j-\frac{1}{2}}} \prod_{h\in H} \frac{\mathcal{X}(z_h)}{q^{z_h-\frac{1}{2}}} \notag \\
\sum_{\substack{n,m \\ Nn = Mm \\ (MNmn,q)=1}} \tau(S\setminus S_J\cup Z_H^- ; n) \tau(Z\setminus Z_H\cup S_J^- ;m),\label{eqn: unitarymoment2}
\end{align}
where, for a multiset $R=\{r_1,\dots,r_k\}$ and a positive integer $n$, $\tau(R;n)$ is defined by
\begin{equation*}
\tau(R;n) := \sum_{n_1\cdots n_k=n} n_1^{-r_1} \cdots n_k^{-r_k}.
\end{equation*}

To evaluate \eqref{eqn: unitarymoment2} further, observe that the orthogonality of Dirichlet characters and M\"{o}bius inversion implies (see also Lemma~2 of \cite{cis2019})
\begin{equation*}
\sumflat_{\chi\mod q}1 = \frac{1}{2}\sum_{d|q} \mu(d)\phi\left( \frac{q}{d}\right) +O(1).
\end{equation*}
Thus, after some rearranging, the right-hand side of \eqref{eqn: unitarymoment2} is, up to an error term that we ignore,
\begin{align}
\frac{1}{2} \sum_{\ell=0}^k \sum_{\substack{J,H \subseteq \{1,\dots,k\} \\ |J|=|H|=\ell }} \prod_{j\in J}\mathcal{X}(s_j) \prod_{h\in H} \mathcal{X}(z_h) \sum_{\substack{n,m \\ Nn = Mm }} \tau(S\setminus S_J\cup Z_H^- ; n) \tau(Z\setminus Z_H\cup S_J^- ;m)\notag\\
\times \sum_{\substack{q\geq 1 \\ (MNmn,q)=1}}g\bfrac qQ \sum_{d|q} \mu(d)\phi\left( \frac{q}{d}\right) q^{\ell -\sum_{j\in J} s_j - \sum_{h\in H} z_h }. \label{U_k}
\end{align}
By Mellin inversion and multiplicativity, we have for any integer $r$ and any complex number $z$ that
\begin{align*}
& \sum_{\substack{q\geq 1 \\ (r,q)=1}}g\bfrac qQ \sum_{d|q} \mu(d)\phi\left( \frac{q}{d}\right) q^{-z }\\
& = \frac{1}{2\pi i } \int_{( 3+ |\text{Re}(z)|)} Q^w \tilde{g} (w)  \frac{\zeta(w+z-1) }{\zeta^2(w+z)}\prod_{p|r}\left(1- \frac{1}{p^{w+z-1}} \right) \left(1- \frac{1}{p^{w+z}}\right)^{-2} \,dw.
\end{align*}
Moving the line of integration to the left, we see from the residue theorem that
\begin{align*}
\sum_{\substack{q\geq 1 \\ (r,q)=1}}g\bfrac qQ \sum_{d|q} \mu(d)\phi\left( \frac{q}{d}\right) q^{-z} \sim   \frac{Q^{2-z} \tilde{g} (2-z)}{\zeta^2(2)}\prod_{p|r}\left(1- \frac{1}{p} \right) \left(1- \frac{1}{p^{2}}\right)^{-2}.
\end{align*}
This, \eqref{eqn: unitarymoment2}, and \eqref{U_k} lead to the following prediction.

\begin{conjecture}[CFKRS~\cite{cfkrs}]\label{Recipe conjecture unitary} Let $\mathcal{U}_{k}(Q;S,Z;M,N)$ be as defined in \eqref{eqn: unitarymoment}. If $\lab\re (s)-\frac{1}{2}\rab\ll 1/\log Q$ and $\lab\re (z)-\frac{1}{2}\rab\ll 1/\log Q$ for all $s\in S$ and $z\in Z$, then as $Q\rightarrow \infty$ we have
\begin{align}
&\mathcal{U}_{k}(Q;S,Z;M,N) \notag\\
&\sim\sum_{\ell=0}^k \sum_{\substack{J,H \subseteq \{1,\dots,k\} \\ |J|=|H|=\ell }} Q^{2+\ell-\sum\limits_{j\in J}s_j-\sum\limits_{h\in H}z_h} \tilde{g}\lz 2+\ell-\sum_{j\in J}s_j-\sum_{h\in H}z_h\pz\cdot \prod_{j\in J}\mathcal{X}(s_j) \prod_{h \in H}\mathcal{X}(z_h) \notag\\
&\hspace{250pt}\times B_{M,N}(S\setminus S_J\cup Z_H^-; Z\setminus Z_H\cup S_J^-), \label{eqn: recipe prediction unitary}
\end{align}
where $S_J=\{s_j : j\in J\}$, $S_J^{-}=\{1-s_j:j\in J\}$, $Z_H=\{z_h : h\in H\}$, $Z_H^{-}=\{1-z_h:h\in H\}$, and
\begin{equation}\label{eqn: Bdef}
B_{M,N}(S,Z):=\frac{1}{2\zeta^2(2)}  \sum_{\substack{n_1,\dots,n_{k},m_1,\dots,m_{k}\geq 1 \\ Nn_1\cdots n_{k}= Mm_{1}\cdots m_{k} }} \frac{1}{n_1^{s_1}\cdots n_{k}^{s_{k}} m_1^{z_1}\cdots m_{k}^{z_{k}} } \prod_{p|MNn_1\cdots m_{k}} \left(1-\frac{1}{p^2} \right)^{-2}\left(1-\frac{1}{p} \right)
\end{equation}
for $S=\{s_1,\dots,s_k\}$ and $Z=\{z_1,\dots,z_k\}$.
\end{conjecture}

\begin{remark}\label{remark: unitary recipe}
\begin{enumerate}
    \item The sum defining $B_{M,N}(S,Z)$ is not absolutely convergent in the proximity of the central point $s_j=z_j=1/2$. It should be interpreted as its meromorphic continuation, which can be obtained in the usual fashion by writing it as an Euler product (see, for example, (2.2) and (4.3) of \cite{btb}). Even though each term on the right-hand side of \eqref{eqn: recipe prediction unitary} has poles, the sum of all the terms is holomorphic in the region with $\lab\re(s)-\frac{1}{2}\rab\ll 1/\log Q$ and $\lab\re(z)-\frac{1}{2}\rab\ll 1/\log Q$ for all $s\in S$ and $z\in Z$. This follows from an application of Lemma~2.5.3 of CFKRS~\cite{cfkrs}.
    \item In the conjecture, we consider the imaginary parts of the variables $s_j,z_j$ as fixed, even though it is believed the conjectures hold with some uniformity in the imaginary parts too. 
    \item Note that $B_{M,N}(S,Z)$ does not depend on the order of the $s_j$'s nor that of the $z_j$'s.
    \item As in \cite{cfkrs}, we expect a power savings error term in the above conjecture.
\end{enumerate} 
\end{remark}

\subsection{Unitary family: predictions via multiple Dirichlet series}\label{sec: unitary MDS}
This discussion easily generalizes to the situation where the number of $s_j$ variables is not the same as the number of $z_j$ variables. However, we assume that there are as many $s_j$ variables as $z_j$ variables for convenience. We apply Mellin inversion to deduce from \eqref{eqn: unitarymoment} that
\begin{equation}\label{U_k as integral}
\mathcal{U}_k(Q;S,Z;M,N) =\frac{1}{2\pi i}\int_{(2+\varepsilon)}A_{M,N}(s_1,\dots,s_k;z_1,\dots,z_k;w)\tilde g(w)Q^w\,dw,
\end{equation} 
where
\begin{equation}\label{eqn: Adef}
\begin{aligned}
A_{M,N}&(s_1,\dots,s_k;z_1,\dots,z_k;w)\\
&:=\sum_{q\geq 1}\sumflat_{\chi\mod q}\frac{\overline{\chi}(M) \chi (N)L(s_1,\chi)\dots L(s_k,\chi)L(z_1,\overline\chi)\dots L(z_k,\overline\chi)}{q^w}.
\end{aligned}
\end{equation}

We first prove some basic analytic properties of this multiple Dirichlet series.

\begin{prop}\label{prop: Aproperties} 
    The multiple Dirichlet series \eqref{eqn: Adef} is absolutely convergent if $\re(w)>2$ and $\re(s_j)>1$ and $ \re(z_j)>1$ for all $j$. It has a meromorphic continuation to the region $\re(w)>1$, $\re(s_j)>2,\ \re(z_j)>2$, where it has a simple pole at $w=2$ with residue 
    $$
    \res{w=2}A_{M,N}(s_1,\dots,s_k;z_1,\dots,z_k;w)=B_{M,N}(S,Z),
    $$
    where $B_{M,N}(S,Z)$ is defined by \eqref{eqn: Bdef}.
\end{prop}
\begin{proof}

The absolute convergence for $\re(w)>2$, $\re(s_j)>1$, $\re(z_j)>1$ follows immediately from the absolute convergence of the Dirichlet series of $L(s,\chi)$ for $\re(s)>1$ and the fact that $\sum^{\flat}_{\chi \bmod{q}}1 \ll q$.
 
We expand each $L(s,\chi)$ in \eqref{eqn: Adef} into its Dirichlet series, and then evaluate the $\chi$-sum via M\"{o}bius inversion and the orthogonality of characters (see Lemma~2 of \cite{cis2019}) to deduce that
\begin{align}
A_{M,N}(s_1,\dots,s_k;z_1,\dots,z_k;w) =\frac{1}{2} \sum_{\substack{q,n_1,\dots,n_{k},m_1,\dots,m_k \\ (q,MNn_1\cdots n_{k}m_1\cdots m_k)=1 }} \frac{1}{q^wn_1^{s_1}\dots n_k^{s_k}m_{1}^{z_{1}}\dots m_k^{z_k}}  \notag\\
\times \sum_{\substack{d|q \\ d| Nn_1\cdots n_k \pm Mm_{1}\cdots m_{k} }} \varphi(d) \mu\left( \frac{q}{d}\right), \label{eqn: applylemma2cis}
\end{align}
where the symbol $\pm$ on the right-hand side means that the right-hand side is a sum of two copies of itself: one with the symbol $\pm$ replaced by $+$, and the other with $\pm$ replaced by $-$. We split \eqref{eqn: applylemma2cis} into the sum of the diagonal terms plus the sum of the off-diagonal terms, i.e., we write it as
\begin{align}
A_{M,N}(s_1,\dots,s_k;z_1,\dots,z_k;w) =A_{D}(s_1,\dots,s_k;z_1,\dots,z_k;w)+A_{OD}(s_1,\dots,s_k;z_1,\dots,z_k;w), \label{eqn: applylemma2cis2}
\end{align}
where
\begin{equation*}
A_D(s_1,\dots,s_k;z_1,\dots,z_k;w)= \frac{1}{2}\sum_{\substack{q,n_1,\dots,n_{k},m_1,\dots,m_k\\ Nn_1\dots n_k=Mm_1\dots m_k \\ (q,Nn_1\cdots n_{k})=1 }} \frac{1}{q^wn_1^{s_1}\dots n_k^{s_k}m_{1}^{z_{1}}\dots m_k^{z_k}}  \sum_{d|q } \varphi(d) \mu\left( \frac{q}{d}\right)
\end{equation*}
and $A_{OD}(s_1,\dots,s_k;z_1,\dots,z_k;w)$ is the sum of the rest of the terms in \eqref{eqn: applylemma2cis}.

To bound the latter sum, observe that by the divisor bound and the fact that $\varphi(d)\leq \varphi(n)$ for $d|n$,
$$
\sum_{\substack{d|q\\d|Nn_1\dots n_k\pm Mm_1\dots m_k}}\phi(d)\mu\bfrac{q}{d}\ll_{\varepsilon} \phi(|Nn_1\dots n_k \pm Mm_1\dots m_k|)(NMn_1\dots n_k m_1\dots m_k)^{\varepsilon}.
$$ 
Since $\phi(|n- m|)\leq |n- m| \leq nm$ for positive integers $n\neq m$, it follows that $A_{OD}$ is absolutely convergent if the variables satisfy $\re(w)>1$, $\re(s_j)>2$, and $\re(z_j)>2$ for all $j$.

Evaluating the sum over $q$ in \eqref{eqn: applylemma2cis2} using multiplicativity, we see that $A_{D}$  equals
\begin{align*}
\frac{\zeta(w-1)}{2\zeta(w)^2}\sum_{\substack{n_1,\dots,n_{k},m_1,\dots,m_k\\ Nn_1\dots n_k=Mm_1\dots m_k}}\frac{1}{n_1^{s_1}\dots m_k^{z_k}}\prod_{p|Nn_1\dots n_k}\lz 1-\frac{1}{p^w}\pz^{-2}\lz1-\frac{1}{p^{w-1}}\pz
\end{align*}
if the variables satisfy $\re(w)>2$, $\re(s_j)>1$, and $\re(z_j)>1$ for all $j$, and it has a meromorphic continuation to $\re(w)>1$, $\re(s_j)>1$, and $\re(z_j)>1$ for all $j$. Thus, $A_{D}$ is a meromorphic function in this region with a simple pole at $w=2$. Its residue at this pole is
\begin{align*}
\lim_{w\rightarrow 2} \frac{(w-2)\zeta(w-1)}{2\zeta(w)^2}\sum_{\substack{n_1,\dots,n_{k},m_1,\dots,m_k\\ Nn_1\dots n_k=Mm_1\dots m_k}}\frac{1}{n_1^{s_1}\dots m_k^{z_k}}\prod_{p|Nn_1\dots n_k}\lz 1-\frac{1}{p^w}\pz^{-2}\lz1-\frac{1}{p^{w-1}}\pz,
\end{align*}
which equals $B_{M,N}(S,Z)$ by \eqref{eqn: Bdef}.
\end{proof}

\begin{remark}\label{remark: extension of residue of MDS}\begin{enumerate}
    \item From Remark \ref{remark: unitary recipe}, we see that the residue $B_{M,N}(S,Z)$ has a meromorphic continuation to a larger region. From the uniqueness of holomorphic extension, it follows that if $(w-2)A_{M,N}(s_1,\dots,s_k;z_1,\dots,z_k;w)$ can be holomorphically extended, then the residue of $A_{M,N}$ at the pole $w=2$ will be given by the meromorphic extension of $B_{M,N}(S,Z)$ (up to possible poles of $B_{M,N}(S,Z)$, i.e., higher order poles of $A(s_1,\dots,s_k;z_1,\dots,z_k;w)$).
    \item We see from the proof that the residue comes from the diagonal terms, while the part that is hindering us from obtaining the meromorphic continuation of $A_{M,N}$ beyond the region with $\re(w)>1$, $\re(s_j)>2$, and $\re(z_j)>2$ are the off-diagonal terms. This is in complete analogy with the recipe.
    \end{enumerate}
\end{remark}

In order to predict an asymptotic formula for \eqref{eqn: unitarymoment} using multiple Dirichlet series, the strategy is to move the line of integration in \eqref{U_k as integral} to the left and collect the residues from the poles. By Proposition~\ref{prop: Aproperties}, the residue at $w=2$ of the integrand in \eqref{U_k as integral} is
$$
Q^2 \tilde{g}(2) B_{M,N}(S,Z).
$$
This is exactly the term corresponding to $\ell=0$ in Conjecture~\ref{Recipe conjecture unitary}. In what follows, we predict that all the other residues come from the pole at $w=2$ shifted by the functional equations of the multiple Dirichlet series.

Given any pair of subsets $J,H\subseteq\{1,\dots,k\}$, using the functional equations \eqref{eqn: fe unitary} of $L(s_j,\chi)$ and $L(z_h,\overline{\chi})$ for each $j\in J$ and $h\in H$ in the definition \eqref{eqn: Adef} gives the functional equation 
$$\begin{aligned}&
A_{M,N}(s_1,\dots,s_k;z_1,\dots,z_k;w)\\
&=\prod_{j\in J}\mathcal{X}(s_j) \prod_{h \in H}\mathcal{X}(z_h) A_{M,N}^{J,H}\lz s_1^J,\dots,s_k^J; z_1^H,\dots, z_k^H;w+\sum_{j\in J}s_j+\sum_{h\in H}z_h-|J|/2-|H|/2\pz,
\end{aligned}$$
where
\begin{align*}
&A_{M,N}^{J,H}\lz s_1,\dots,s_k; z_1,\dots, z_k;w\pz\\
&:=\sum_{q\geq 1}\sideset{}{^\flat}\sum_{\chi\mod q}\frac{\varepsilon(\chi)^{|J|-|H|}\overline{\chi}(M) \chi (N)}{q^w}\prod_{j\notin J}L(s_j,\chi)\prod_{h\in H}L(z_h,\chi)\prod_{j\in J} L(s_j,\overline\chi) \prod_{h\notin H} L(z_h,\overline\chi).
\end{align*}
Let
$$
\sigma_{J,H}(w):=w+\sum_{j\in J}s_j+\sum_{h\in H}z_h-|J|/2-|H|/2.
$$
If $|J|=|H|$, then $A_{M,N}^{J,H}$ is equal to $A_{M,N}$ up to a certain transformation of the variables. Thus, in this case, the corresponding functional equation shifts the pole of $A_{M,N}$ at $w=2$ to a pole at $\sigma_{J,H}(w)=2$ whose residue is
\begin{align*}
& \prod_{j\in J}\mathcal{X}(s_j) \prod_{h \in H}\mathcal{X}(z_h) \ \res{w=2}A_{M,N}^{J,H}\lz s_1^J,\dots,s_k^J; z_1^H,\dots z_k^H;w\pz\\
&= B_{M,N}(S\setminus S_J\cup Z_H^-; Z\setminus Z_H\cup S_J^-) \prod_{j\in J}\mathcal{X}(s_j) \prod_{h \in H}\mathcal{X}(z_h)
\end{align*}
by Proposition~\ref{prop: Aproperties}. Hence,  if we denote $\ell=|J|=|H|$, then the residue of the integrand in \eqref{U_k as integral} at $\sigma_{J,H}(w)=2$ is
\begin{align}
Q^{2+\ell-\sum\limits_{j\in J}s_j-\sum\limits_{h\in H}z_h}
& \tilde g\lz 2+\ell-\sum_{j\in J}s_j-\sum_{h\in H}z_h\pz \notag\\
& \times B_{M,N}(S\setminus S_J\cup Z_H^-; Z\setminus Z_H\cup S_J^-) \prod_{j\in J}\mathcal{X}(s_j) \prod_{h \in H}\mathcal{X}(z_h). \label{eqn: unitaryresidue}
\end{align}
Observe that this is the same as the term corresponding to the pair $J,H$ in Conjecture~\ref{Recipe conjecture unitary}.

On the other hand, if $|J|\neq |H|$, then the terms in the definition of $A_{M,N}^{J,H}$ are twisted by the factor $\varepsilon(\chi)^{|J|-|H|}$. Similarly as in the recipe, we conjecture that the oscillation of the $\varepsilon(\chi)$ makes it so that the $A_{M,N}^{J,H}$ with $|J|\neq |H|$ do not contribute to the main terms with additional poles. 

We thus state our conjecture as follows

\begin{conjecture}\label{MDS conjecture unitary}
    The function
    \begin{equation}\label{eqn: MDS conjecture unitary}
    \Bigg(\prod_{\substack{J,H\subset\{1,\dots,k\}\\|J|=|H|}}\lz\sigma_{J,H}(w)-2\pz \Bigg)\Bigg( A_{M,N}(s_1,\dots,s_k;z_1,\dots,z_k;w) - \prod_{j=1}^k \zeta(s_j)\zeta(z_j)\Bigg)
    \end{equation}
    has a holomorphic continuation to a tube domain that contains the point
    $$
     (s_1,\dots,s_k,z_1,\dots,z_k,w)=(\tfrac{1}{2},\dots,\tfrac{1}{2},2),
    $$ and it is polynomially bounded in vertical strips in this region.
\end{conjecture}
\begin{remark}
\begin{enumerate}
\item Here and onwards, polynomial boundedness of $f(u_1,\dots,u_k)$ in vertical strips means that for any fixed real numbers $a_1,\dots,a_k,b_1,\dots,b_k$, there is a polynomial $P(t_1,\dots,t_k)$, such that
    $$
    |f(\sigma_1+it_1,\dots,\sigma_k+it_k)|\leq |P(t_1,\dots,t_k)|
    $$
    for all real $t_1,\dots,t_k$ and $\sigma_1,\dots,\sigma_k$ such that $a_j\leq \sigma_j\leq b_j$ for all $j$.
    \item The product $\prod_{j=1}^k \zeta(s_j)\zeta(z_j)$ is the $q=1$ term in the definition \eqref{eqn: Adef}. We take the Dirichlet character modulo $1$ to be primitive by convention. We subtract this product of zeta functions from $A_{M,N}$ in \eqref{eqn: MDS conjecture unitary} in order to remove the polar divisors of $A_{M,N}$ at the hyperplanes $s_j=1$ and $z_j=1$.
    \item Conjecture~\ref{MDS conjecture unitary} is analogous to Conjecture~2.7 in \cite{dgh}, since, like Conjecture~2.7 of \cite{dgh}, it asserts that the poles of $A$ are only those that are shifts of the pole at $w=2$ by the functional equation that sends $w$ to $\sigma_{J,H}(w)$ with $|J|=|H|$. Indeed, in \cite{dgh}, Proposition~2.6 together with Conjecture~2.7 assert that the poles of $Z_{\epsilon_1,\dots,\epsilon_{2m},0}$ there are only those that are reflections of $w=1$ under the functional equations given by applying the functional equation of the Riemann zeta-function $\zeta(s)$ to the same number of factors $\zeta(s_j+it)$ as factors $\zeta(s_j-it)$ in the definition of $Z_{\epsilon_1,\dots,\epsilon_{2m},0}$.
    \end{enumerate}
\end{remark}

\begin{prop}
Conjecture~\ref{MDS conjecture unitary} implies Conjecture~\ref{Recipe conjecture unitary} with a power-saving error term.
\end{prop}
\begin{proof}
Let $s_1,\dots,s_k,z_1,\dots,z_k$ be as in the statement of Conjecture~\ref{Recipe conjecture unitary}. Assume that $Q$ is large and $\delta$ is small, so that $(s_1,\dots,s_k;z_1,\dots,z_k;2-\delta)$ is in the tube domain from the statement of Conjecture \ref{MDS conjecture unitary}. We move the line of integration in \eqref{U_k as integral} to $\re(w)=2-\delta$ for some $\delta>0$. This leaves residues from the poles at $\sigma_{J,H}(w)=2$, where $J,H$ range over the subsets of $\{1,\dots,k\}$ such that $|J|=|H|$. By Proposition \ref{prop: Aproperties} and Remark \ref{remark: extension of residue of MDS}, the residues at these poles are given by \eqref{eqn: unitaryresidue} and exactly correspond to the main terms from Conjecture \ref{Recipe conjecture unitary}. The integral along $\re(w)=2-\delta$ is $O(Q^{2-\delta})$ by the rapid decay of $\tilde{g}$ and the polynomial growth of $A_{M,N}$ assumed in Conjecture~\ref{MDS conjecture unitary}.
\end{proof}
\begin{remark}
A similar proof shows that Conjecture \ref{MDS conjecture unitary} implies the asymptotic formula in Conjecture \ref{Recipe conjecture unitary} in a wider range of the variables $\re(s_j),\re(z_j)>1/2-\epsilon$ for some $\epsilon>0$.
\end{remark}

\section{Symplectic family: Dirichlet $L$-functions of real primitive characters}\label{sec: symplectic}

This section also forms part of the second author's thesis \cite{thesis}, and was recently studied in more detail in \cite{cech2}, where explicit regions of meromorphic continuation of the corresponding multiple Dirichlet series are determined. In this section, we write $\chi_d$ for the Kronecker symbol $\bfrac{d}{\cdot}$, where $d$ denotes a positive fundamental discriminant. For each positive fundamental discriminant $d$, $L(s,\chi_d)$ satisfies the functional equation
\begin{equation}\label{eqn: fe symplectic}
	L(s,\chi_d)=d^{1/2-s}\mathcal{X}(s) L(1-s,\chi_d),
\end{equation}
where $\mathcal{X}$ is from \eqref{eqn: Xdef}. We also have the approximate functional equation
\begin{equation} \label{eqn: approximate fe symplectic}
L(s,\chi_d)\approx\sum_{n}\frac{\chi_d(n)}{n^{s}}+d^{\frac{1}{2}-s} \mathcal{X}(s)\sum_{n}\frac{\chi_d(n)}{n^{1-s}}.
\end{equation}

Our aim is to understand the twisted shifted moment
\begin{equation}\label{eqn: symplecticmoment}
	\mathcal{S}_{k}(X;S;M)=\sumstar_{d\geq 1} g\left( \frac{d}{X}\right)\chi_d(M)L(s_1,\chi_d)\dots L(s_k,\chi_d),
\end{equation}
where $X$ is a parameter tending to infinity, $k,M$ are positive integers, and $S=\{s_1,\dots,s_k\}$ is a multiset of complex numbers close to $1/2$, say $s_j-\frac{1}{2}\ll (\log X)^{-1}$ for each $j$. The summation symbol $\sum^*$ here means that we are summing over fundamental discriminants $d$.

\subsection{Symplectic family: carrying out the CFKRS recipe}

We replace each $L(s,\chi_d)$ in \eqref{eqn: symplecticmoment} by its approximate functional equation \eqref{eqn: approximate fe symplectic}. We multiply out the resulting product to obtain the approximation
\begin{equation}\label{eqn: sjj}
\mathcal{S}_{k}(X;S;M)\approx \sum_{J\subseteq \{1,\dots,k\}} \sumstar_{d\geq 1} g\left( \frac{d}{X}\right) \chi_d(M) \prod_{j\in J}\frac{\mathcal{X}(s_j) }{d^{s_j-\frac{1}{2}}} \sum_{n_1,\dots,n_k}\frac{\chi_d(n_1\cdots n_k)}{n_1^{s_1^J}\dots n_k^{s_k^J}},
\end{equation}
where $s_j^J$ is defined by \eqref{eqn: sjjdef}. 

By \cite[(3.1.21)]{cfkrs} and the fact that the number of positive fundamental discriminants $\leq y$ is asymptotic to $y/(2\zeta(2))$, we have
\begin{equation*}
	\sumstar_{1\leq d\leq y}\chi_d(M n_1\dots n_k)\approx\begin{cases}
	    \displaystyle \frac{y}{2\zeta(2)}a(M n_1\dots n_k),&\hbox{if $M n_1\dots n_k=\square$,}\\
     0,&\hbox{if $M n_1\dots n_k\neq\square$,}
	\end{cases}
\end{equation*}
where
\begin{equation}\label{eqn: a(n)def}
a(n):=\prod_{p|n}\frac{p}{p+1}.
\end{equation}
From this, Mellin inversion, and partial summation, we expect for any complex number $z$ that
\begin{equation}\label{eqn: dsumaverage}
\sumstar_{d\geq 1} g\left( \frac{d}{X}\right) \frac{\chi_d(Mn_1\cdots n_k) }{d^z}
\approx \frac{a(Mn_1\cdots n_k)\mathbf{1}_{Mn_1\cdots n_k=\square}}{4\pi i\zeta(2)}\int_{(2+|\text{Re}(z)|)}  \frac{X^w\tilde{g}(w)}{w+z-1}  \,dw,
\end{equation}
where $\mathbf{1}_{Mn_1\cdots n_k=\square}=1$ if $Mn_1\cdots n_k$ is a square and is zero otherwise. Moving the line of integration to the left and using the residue theorem, we arrive at the approximation
$$
\sumstar_{d\geq 1} g\left( \frac{d}{X}\right) \frac{\chi_d(Mn_1\cdots n_k) }{d^z}   \approx  \frac{a(Mn_1\cdots n_k) \mathbf{1}_{Mn_1\cdots n_k=\square}}{2\zeta(2)} X^{1-z} \tilde{g}(1-z).
$$
This and \eqref{eqn: sjj} lead to the following prediction.
\begin{conjecture}[CFKRS~\cite{cfkrs}]\label{Recipe conjecture symplectic} Let $\mathcal{S}_{k}(X;S;M)$ be as defined in \eqref{eqn: symplecticmoment}. If $\lab\re(s)-\frac{1}{2}\rab\ll 1/\log X$ for all $s\in S$, then as $X\rightarrow \infty$ we have
\begin{align}\label{eqn: recipe prediction symplectic}
\mathcal{S}_{k}(X;S;M)\sim \sum_{J\subseteq\{1,\dots,k\}}\prod_{j\in J}\mathcal{X}(s_j)\cdot X^{1+\frac{|J|}{2}-\sum\limits_{j\in J}s_j}\tilde g\lz1+\frac{|J|}{2}-\sum_{j\in J}s_j\pz T_M(S\setminus S_J\cup S_J^-),
\end{align}
where $S_J=\{s_j:j\in J\}$, $S_J^{-}=\{1-s_j :j\in J\}$, and
\begin{equation}\label{eqn: TMdef}
T_M(S):=\frac{1}{2\zeta(2)}\sum_{M n_1\dots n_k=\square}\frac{a(M n_1\dots n_k)}{n_1^{s_1}\dots n_k^{s_k}}
\end{equation}
for $S=\{s_1,\dots,s_k\}$, with $a(n)$ defined by \eqref{eqn: a(n)def}.
\end{conjecture}

\subsection{Symplectic family: predictions via multiple Dirichlet series}
We apply Mellin inversion to deduce from \eqref{eqn: symplecticmoment} that
\begin{equation}\label{Moments after Perron}
\mathcal{S}_k(X;S;M)=\frac{1}{2\pi i}\int_{(2)}A_M(s_1,\dots,s_k;w)X^w\tilde g(w)dw,
\end{equation}
where
\begin{equation}\label{eqn: symplecticAdef}	A_M(s_1,\dots,s_k;w):=\sumstar_{d\geq 1}\frac{\chi_d(M)L(s_1,\chi_d)\dots L(s_k,\chi_d)}{d^w}.
\end{equation}

We first prove some basic analytic properties of this multiple Dirichlet series (this is similar to Propositions 3.2 and 3.3 of \cite{dgh}).

\begin{prop}\label{prop: symplectic region}
    The multiple Dirichlet series \eqref{eqn: symplecticAdef} is absolutely convergent if $\re(w)>1$ and $\re(s_j)>1$ for all $j$. It has a meromorphic continuation to the region
    $\re(w)>1/2$, $\re(s_j)>5/4$, where it has a simple pole at $w=1$ with residue
    $$
    \res{w=1}A_M(s_1,\dots,s_k;w)=T_M(S),$$
    where $T_M(S)$ is defined by \eqref{eqn: TMdef}.
\end{prop}
\begin{proof}
The absolute convergence for $\re(w)>1$ and $\re(s_j)>1$ follows immediately from the absolute convergence of the Dirichlet series of $L(s,\chi)$ for $\re(s)>1$. 

We expand each $L(s,\chi_d)$ in \eqref{eqn: symplecticAdef} into its Dirichlet series to deduce that
\begin{equation}\label{Z exchanging sums}	A_M(s_1,\dots,s_k;w)
=\sum_{n_1,\dots,n_k\geq 1}\frac{L_D\lz w,\bfrac{\cdot}{M n_1\dots n_k}\pz}{n_1^{s_1}\dots n_k^{s_k}},
\end{equation}
where
\begin{equation*}
	L_D(w,\chi):=\sumstar_{d\geq 1}\frac{\chi(d)}{d^w}.
\end{equation*} 
By \cite[Lemma~2.4]{cech} and the convexity bound, we have for $\re(w)>1/2$, away from the possible pole at $w=1$,
$$
L_D\lz w,\bfrac{\cdot}{M n_1\dots n_k}\pz \ll (|w| Mn_1\dots n_k)^{1/4+\epsilon}.
$$
This and \eqref{Z exchanging sums} imply that \eqref{eqn: symplecticAdef} has a meromorphic continuation to the region $\re(w)>1/2$, $\re(s_j)>5/4$. In this region, the only possible polar divisor is at $w=1$. Indeed, by Lemma~2.4 and equations (4.8) and (4.9) of \cite{cech}, if $n$ is not a square, then $L_D\lz w,\bfrac{\cdot}{n}\pz$ is holomorphic in the region $\re(w)>1/2$, while if $n$ is a square, then it is holomorphic there except for a simple pole at $w=1$. The residue at this pole is $a(n)/(2\zeta(2))$. Hence, for any fixed $s_1,\dots,s_k$ with real parts $>5/4$, \eqref{Z exchanging sums} has a pole at $w=1$ with residue $T_M(S)$.
\end{proof}

By Proposition \ref{prop: symplectic region}, we see that the integrand in \eqref{Moments after Perron} has a simple pole at $w=1$ with residue
$$X\tilde g(1) T_M(S),$$ which exactly corresponds to the diagonal term $J=\emptyset$ in the recipe prediction \eqref{eqn: recipe prediction symplectic}. We will see that the other terms come from this one shifted by the functional equations.

 Given any subset $J\subseteq \{1,\dots,k\}$, applying the functional equation \eqref{eqn: fe symplectic} to $L(s_j,\chi_d)$ for each $j\in J$ in \eqref{eqn: symplecticAdef} gives
\begin{equation}\label{eqn: symplecticfunceqn}
A_M(s_1,\dots,s_k;w)= \prod_{j\in J} \mathcal{X}(s_j) \ A_M\left(s_1^J,\dots,s_k^J; w + \sum_{j\in J} s_j -\frac{|J|}{2} \right), 
\end{equation}
where $s_j^J$ is defined by \eqref{eqn: sjjdef}. 

Let $$\sigma_J(w):=w+\sum_{j\in J}s_j-\frac{|J|}{2}.$$
The functional equation \eqref{eqn: symplecticfunceqn} shifts the pole at $w=1$ to a pole at 
\begin{equation*}
	\sigma_J(w)=1,
\end{equation*}
whose residue is
\begin{equation}\label{eqn: symplecticresidue}
\prod_{j\in J} \mathcal{X}(s_j) \ \res{w=1} A_M(s_1^J,\dots,s_k^J; w) =\prod_{j\in J} \mathcal{X}(s_j)\cdot T_M(S\setminus S_J\cup S_J^{-}) 
\end{equation}
by Proposition~\ref{prop: symplectic region}.

Thus, the residue of the integrand in \eqref{Moments after Perron} at this pole is
\begin{equation*}
X^{1+\frac{|J|}{2}-\sum_{j\in J}s_j} \tilde{g} \left( 1+\frac{|J|}{2}-\sum_{j\in J}s_j \right) T_M(S\setminus S_J\cup S_J^{-}) \prod_{j\in J} \mathcal{X}(s_j),
\end{equation*}
which is the same as the term corresponding to $J$ in Conjecture~\ref{Recipe conjecture symplectic}. We conjecture that these poles are the only poles close to the point $(\frac{1}{2},\dots,\frac{1}{2},1)$.

\begin{conjecture}\label{MDS conjecture symplectic}
    The function
    $$\Bigg(\prod_{J\subseteq\{1,\dots,k\}}(\sigma_J(w)-1) \Bigg)\Bigg( A_M(s_1,\dots,s_k;w) -\prod_{j=1}^k \zeta(s_j)\Bigg)$$
    has a holomorphic continuation to a tube domain that contains the point
    $$
    (s_1,\dots,s_k,w) = (\tfrac{1}{2},\dots,\tfrac{1}{2},1),
    $$
    and it is polynomially bounded in vertical strips in this region.
\end{conjecture}

\begin{prop}\label{prop: MDS implies recipe symplectic}
Conjecture~\ref{MDS conjecture symplectic} implies Conjecture~\ref{Recipe conjecture symplectic} with a power-saving error term.
\end{prop}
\begin{proof}
Let $s_1,\dots,s_k$ be as in the statement of Conjecture~\ref{Recipe conjecture symplectic}. We move the line of integration in \eqref{Moments after Perron} to $\re(w)=1-\delta$ for some $\delta>0$. This leaves residues from the poles at the points $\sigma_J(w)=1$ for $J\subseteq\{1,\dots,k\}$. As we have shown above, the residues at these poles are given by \eqref{eqn: symplecticresidue} and exactly correspond to the terms from Conjecture \ref{Recipe conjecture symplectic}. The integral along $\re(w)=1-\delta$ is $O(X^{1-\delta})$ by the rapid decay of $\tilde{g}$ and the polynomial growth of $A$ assumed in Conjecture~\ref{MDS conjecture symplectic}.
\end{proof}

Let us now compare Conjecture~\ref{MDS conjecture symplectic} and Proposition~\ref{prop: MDS implies recipe symplectic} above with Conjecture 3.6 and Theorem 3.7 of \cite{dgh}. Diaconu, Goldfeld and Hoffstein~\cite{dgh} use a (unpublished) version of a Tauberian-type theorem due to Stark. This enables them to find the main term in the asymptotic formula for the $m$th (unshifted) moment at the central point assuming that the corresponding multiple Dirichlet series has a meromorphic continuation. No polynomial boundedness hypothesis is needed to deduce the asymptotic formula through this method (see Conjecture~3.6 of \cite{dgh}). In our case, we consider the shifted moments and  apply Mellin inversion and the residue theorem to obtain a direct comparison between the arising residues with all the terms from the CFKRS recipe. As in CFKRS~\cite{cfkrs}, this leads to the prediction of the main term and lower order terms in the asymptotic formula for the unshifted moment, up to a power-saving error term. The polynomial boundedness hypothesis in Conjecture~\ref{MDS conjecture symplectic} facilitates our method.

We should remark that the necessary work establishing the Keating-Snaith conjectures for the unshifted moment at the central point assuming Conjecture~\ref{Recipe conjecture symplectic} is done by CFKRS~\cite{cfkrs} via a limit calculation. Part of the proof of Theorem 3.7 of \cite{dgh} is also a calculation of a limit that reaches the same prediction as Keating and Snaith.

\section{Orthogonal family: quadratic twists of modular Hecke eigenforms}\label{sec: orthogonal}

In this section, we let $f$ be a Hecke eigencuspform for the full modular group $\mathrm{SL}(2,\Z)$ of (even) weight $\kappa$, with Fourier expansion
\begin{equation}\label{eqn: fourier}
f(z)=\sum_{n=1}^\infty a_f(n) n^{\frac{\kappa-1}{2}}e^{2\pi i nz},
\end{equation}
normalized such that $|a_f(n)|\leq \tau(n)$, the divisor function. Recall from the previous section that if $d$ is a fundamental discriminant then we define $\chi_d=\left(\frac{d}{\cdot}\right)$, the Kronecker symbol. For a fundamental discriminant $d$, the twisted $L$-function $L(s,f\otimes \chi_d)$ is defined for $\re(s)>1$ by
\begin{equation*}
L(s,f\otimes\chi_d)=\sum_{n=1}^\infty\frac{a_f(n)\chi_d(n)}{n^s}.
\end{equation*}

For each positive fundamental discriminant $d$, $L(s,f\otimes \chi_d)$ satisfies the functional equation (see, for example, \cite[p.~1098]{soundyoung}) 
\begin{equation}\label{eqn: fe orthogonal}
L(s,f\otimes\chi_d)=d^{1-2s}\mathcal{X}_f(s) L(1-s,f\otimes\chi_d),
\end{equation}
where $\mathcal{X}_f(s)$ is
\begin{equation*}
\mathcal{X}_f(s)=i^{\kappa} (2\pi)^{2s-1}\frac{\Gamma\lz1-s+\frac{\kappa-1}{2}\pz}{\Gamma\lz s+\frac{\kappa-1}{2}\pz}.
\end{equation*}
We also have the approximate functional equation
\begin{equation}\label{eqn: approximate fe orthogonal}
L(s,f\otimes\chi_d) \approx \sum_{n}\frac{a_f(n)\chi_d(n)}{n^s}+d^{1-2s}\mathcal{X}_f(s)\sum_{n}\frac{a_f(n)\chi_d(n)}{n^{1-s}}.
\end{equation}

We wish to understand the twisted shifted moment
\begin{equation}\label{eqn: orthogonalmoment}
\mathcal{O}_{f,k}(X;S;M)=\sumstar_{d\geq 1}g\bfrac{d}{X}\chi_d(M)L(s_1,f\otimes\chi_d)\dots L(s_k,f\otimes\chi_d),
\end{equation}
where $X$ is a parameter tending to infinity, $k,M$ are positive integers, and $S=\{s_1,\dots,s_k\}$ is a multiset of complex numbers close to $1/2$, say $s_j-\frac{1}{2}\ll (\log X)^{-1}$ for each $j$. As in the previous section, the summation symbol $\sum^*$ here means that we are summing over fundamental discriminants $d$.

\subsection{Orthogonal family: carrying out the CFKRS recipe}
To carry out the recipe from \cite{cfkrs}, we replace each $L(s,f\otimes \chi_d)$ in \eqref{eqn: orthogonalmoment} by its approximate functional equation \eqref{eqn: approximate fe orthogonal}. We multiply out the resulting product to obtain the approximation
\begin{equation*}
\begin{aligned}
\mathcal{O}_{f,k}(X;S;M)\approx \sum_{J\subseteq \{1,\dots,k\}}& \sumstar_{d\geq 1} g\left( \frac{d}{X}\right) \chi_d(M) \prod_{j\in J}\frac{\mathcal{X}_f(s_j) }{d^{2s_j-1}}\\
&\cdot\quad\sum_{n_1,\dots,n_k\geq 1}\frac{a_f(n_1)\cdots a_f(n_k)\chi_d(n_1\cdots n_k)}{n_1^{s_1^J}\dots n_k^{s_k^J}},
\end{aligned}
\end{equation*}
where, as before, $s_j^J$ is defined by \eqref{eqn: sjjdef}. From this and a procedure similar to our estimation of \eqref{eqn: sjj}, we are led to the following prediction.
\begin{conjecture}[CFKRS~\cite{cfkrs}]\label{Recipe conjecture orthogonal} Let $\mathcal{O}_{f,k}(X;S;M)$ be as defined in \eqref{eqn: orthogonalmoment}. If $\lab\re(s)-\frac{1}{2}\rab\ll 1/\log X$ for all $s\in S$, then as $X\rightarrow \infty$ we have
\begin{align*}
\mathcal{O}_{f,k}(X;S;M)\sim \sum_{J\subseteq\{1,\dots,k\}}\prod_{j\in J}\mathcal{X}_f(s_j)X^{1+|J|-2\sum\limits_{j\in J}s_j}\tilde g\lz1+|J|-2\sum_{j\in J}s_j\pz H_{f,M}(S\setminus S_J\cup S_J^-),
\end{align*}
where $S_J=\{s_j:j\in J\}$, $S_J^{-}=\{1-s_j :j\in J\}$, and
\begin{equation}\label{eqn: HfMdef}
H_{f,M}(S):=\frac{1}{2\zeta(2)}\sum_{M n_1\dots n_k=\square}\frac{a_f(n_1)\cdots a_f(n_k)a(M n_1\dots n_k)}{n_1^{s_1}\dots n_k^{s_k}}
\end{equation}
for $S=\{s_1,\dots,s_k\}$, with $a_f(n)$ defined by \eqref{eqn: fourier} and $a(n)$ defined by \eqref{eqn: a(n)def}.
\end{conjecture}

\subsection{Orthogonal family: predictions via multiple Dirichlet series}

We apply Mellin inversion to deduce from \eqref{eqn: orthogonalmoment} that
\begin{equation}\label{Moments after Perron2}
\mathcal{O}_{f,k}(X;S;M)=\frac{1}{2\pi i}\int_{(2)}A_{f,M}(s_1,\dots,s_k;w)X^w\tilde g(w)dw,
\end{equation}
where
\begin{equation}\label{eqn: orthogonalAdef}	A_{f,M}(s_1,\dots,s_k;w):=\sumstar_{d\geq 1}\frac{\chi_d(M)L(s_1,f\otimes \chi_d)\dots L(s_k,f\otimes \chi_d)}{d^w}.
\end{equation}

We first prove some basic analytic properties of this multiple Dirichlet series.

\begin{prop}\label{prop: orthogonal region}
    The multiple Dirichlet series \eqref{eqn: orthogonalAdef} is absolutely convergent if $\re(w)>1$ and $\re(s_j)>1$ for all $j$. It has a meromorphic continuation to the region
    $\re(w)>1/2$, $\re(s_j)>5/4$, where it has a simple pole at $w=1$ with residue
    $$
    \res{w=1}A_{f,M}(s_1,\dots,s_k;w)=H_{f,M}(S),
    $$
    where $H_{f,M}(S)$ is defined by \eqref{eqn: HfMdef}.
\end{prop}
\begin{proof}
The proof is similar to that of Proposition~\ref{prop: symplectic region}. Recall from \eqref{eqn: fourier} that $|a_f(n)|\leq \tau(n)$ for all $n$.
\end{proof}

Given any subset $J\subseteq \{1,\dots,k\}$, applying the functional equation \eqref{eqn: fe orthogonal} of $L(s_j,f\otimes \chi_d)$ for each $j\in J$ in \eqref{eqn: orthogonalAdef} gives
\begin{equation}\label{eqn: AfMfunceqn}
A_{f,M}(s_1,\dots,s_k;w)=\prod_{j\in J}\mathcal{X}_{f}(s_j)\cdot A_{f,M} \lz s_1^J,\dots,s_k^J;w+2\sum_{j\in J}s_j-|J|\pz,
\end{equation}
where $s_j^J$ is defined by \eqref{eqn: sjjdef}. Let
$$
\sigma_{J}(w):=w+2\sum_{j\in J}s_j-|J|.
$$
The functional equation \eqref{eqn: AfMfunceqn} shifts the pole at $w=1$ to a pole at
\begin{equation*}
\sigma_J(w)=1,
\end{equation*}
whose residue is
\begin{equation*}
\prod_{j\in J} \mathcal{X}_f (s_j)\  \res{w=1}A_{f,M}(s_1^J,\dots,s_k^J;w)=H_{f,M}(S\setminus S_J\cup S_J^-)\prod_{j\in J} \mathcal{X}_f (s_j)
\end{equation*}
by Proposition~\ref{prop: orthogonal region}.

Thus, the residue of the integrand in \eqref{Moments after Perron2} at this pole is
\begin{equation*}
X^{1+|J|-2\sum\limits_{j\in J}s_j} \tilde{g}\left( 1+|J|-2\sum_{j\in J}s_j\right) H_{f,M}(S\setminus S_J\cup S_J^-)\prod_{j\in J} \mathcal{X}_f (s_j),
\end{equation*}
which is the same as the term corresponding to $J$ in Conjecture~\ref{Recipe conjecture orthogonal}. We conjecture that these are the only poles close to the point $(\frac{1}{2},\dots,\frac{1}{2},1)$.

\begin{conjecture}\label{MDS conjecture orthogonal}
    The function $$\prod_{J\subset\{1,\dots,k\}}(\sigma_J(w)-1)\cdot A_{f,M}(s_1,\dots,s_k;w)$$ has a holomorphic continuation to a tube domain that contains the point
    $$
    (s_1,\dots,s_k,w) = (\tfrac{1}{2},\dots,\tfrac{1}{2},1),
    $$
    and it is polynomially bounded in this region.
\end{conjecture}

\begin{prop}
Conjecture~\ref{MDS conjecture orthogonal} implies Conjecture~\ref{Recipe conjecture orthogonal} with a power-saving error term.
\end{prop}
\begin{proof}
The proof is similar to that of Proposition~\ref{prop: MDS implies recipe symplectic}. Move the line of integration in \eqref{Moments after Perron2} to $\re(w)=1-\varepsilon$.
\end{proof}

\section{Quadratic twists of elliptic curves}\label{sec: elliptic}

The discussion in the previous section can be extended to congruence subgroups, with a complication stemming from the root number. In this section, we present the extension in the special case of $L$-functions associated to elliptic curves over $\mathbb{Q}$.

We fix an elliptic curve $E$ over $\mathbb{Q}$ with conductor $N_E=N$ that is not a square, and define the (normalized) Hasse-Weil $L$-function
$$L(s,E)=\prod_{p\nmid N}\lz1-\frac{\lambda_E(p)}{p^{s}}+\frac{1}{p^{2s}}\pz^{-1}\prod_{p|N}\lz 1-\frac{\lambda_E(p)}{p^{s}}\pz^{-1}=\sum_{n\geq 1}\frac{\lambda_E(n)}{n^{s}},$$ where
$$|\lambda_E(p)|=\frac1{\sqrt p}|p+1-\#E(\mathbb{F}_p)|<2$$ by the Hasse bound.

By the modularity theorem due to Wiles~\cite{wiles}, Taylor and Wiles~\cite{taylorwiles}, and Breuil, Conrad, Diamond, and Taylor~\cite{bcdt}, there is a primitive cuspidal modular newform of weight 2 and level $N$ such that $L(s,f)=L(s,E)$. Therefore $L(s,E)$ has a meromorphic continuation and satisfies the functional equation (see, for example, \cite{radziwillsound})
$$
L(s,E)=\mathcal{X}_E(s)L(1-s,E),
$$
where
$$
\mathcal{X}_E(s)=\varepsilon(E)\bfrac{\sqrt N}{2\pi}^{1-2s}\frac{\Gamma\lz\frac32-s\pz}{\Gamma\lz\frac12+s\pz},
$$
with $\varepsilon(E)=\pm1$ being the root number of $E$.

We will consider the $L$-functions of $E$ twisted by quadratic characters $\chi_d$. For a fundamental discriminant $d$, we define
$$
L(s,E_d)=\sum_{n\geq 1}\frac{\lambda_E(n)\chi_d(n)}{n^{s}}.
$$
If $(d,N)=1$, then we have the functional equation (see, for example, \cite{radziwillsound}) \begin{equation}\label{eqn: elliptic twist functional equation}
L(s,E_d)=|d|^{1-2s}\chi_d(-N)\mathcal{X}_E(s)L(1-s,E_d).
\end{equation}
Here, the root number of $L(s,E_d)$ is $\varepsilon(E)\chi_d(-N)$. Let us note that the correlation between the root numbers $\varepsilon(E)\chi_d(N)$ and the coefficients $\lambda_E(n)\chi_d(n)$ leads to the disagreement between the recipe prediction and the multiple Dirichlet series prediction.

We also have the approximate functional equation
\begin{equation}\label{eqn: approximate fe elliptic}
L(s,E_d)\approx \sum_{n}\frac{\lambda_E(n)\chi_d(n)}{n^s}+|d|^{1-2s}\chi_d(-N)\mathcal{X}_E(s)\sum_n\frac{\lambda_E(n)\chi_d(n)}{n^{1-s}},
\end{equation}

Our goal is to understand the moments
\begin{equation}\label{eqn: ellipticmoment}
\mathcal{O}_{E,k}(X;S;M)=\sumstar_{\substack{d\geq 1 \\ (d,N)=1}}\chi_d(M)L(s_1,E_d)\dots L(s_k,E_d)g\bfrac dX,
\end{equation}
where the sum runs over positive fundamental discriminants $d$ that are coprime to $N$.

\subsection{Quadratic twists of elliptic curves: carrying out the CFKRS recipe} In this section, we apply two different versions of the recipe to predict two different asymptotic formulas for \eqref{eqn: ellipticmoment}. The first version is the exact one outlined in Section~4.1 of CFKRS~\cite{cfkrs}. The prediction of this version is different from the multiple Dirichlet series prediction in the next section (Proposition~\ref{prop: MDSelliptic} below). The second version is a slight modification of the recipe outlined in Section~4.1 of CFKRS~\cite{cfkrs}. In this second version, we multiply the root number and the Dirichlet series coefficients of $\prod_{j=1}^k L(s_j,E_d)$, and then take the average of this product over $d$ instead of estimating the average of the root number and the average of the product of the Dirichlet series coefficients separately. The prediction of this second version of the recipe agrees with the multiple Dirichlet series prediction.

In both versions, we start by replacing each $L$-function in \eqref{eqn: ellipticmoment} with the right-hand side of \eqref{eqn: approximate fe elliptic}
and write
\begin{align}
\mathcal{O}_{E,k}(X;S;M)\approx\prod_{J\subset\{1,\dots,k\}}\sumstar_{\substack{d\geq 1 \\ (d,N)=1}}g\bfrac dX\chi_d(M)\prod_{j\in J}\frac{\mathcal{X}_E(s_j)\chi_d(N)}{d^{2s_j-1}} \notag\\
\cdot\quad \sum_{n_1,\dots,n_k\geq 1}\frac{\lambda_E(n_1)\dots \lambda_E(n_k)\chi_d(n_1\dots n_k)}{n_1^{s_1}\dots n_k^{s_k}}.\label{eqn: elliptic moment after approx fe}
\end{align}
According to Step (3) of the recipe \cite[Section 4.1]{cfkrs}, we next replace the contribution of the root numbers by their expected value as averaged over the family. For any $J\subseteq\{1,\dots,k\}$, this expected value is given by (recall that we assume that $N$ is not a square)
$$\lim\limits_{X\rightarrow\infty}\frac{\varepsilon(E)^{|J|}}{\#\{1\leq d\leq X:\ d\text{ fund. discr., }(d,N)=1\}}\sumstar_{\substack{1\leq d\leq X\\(d,N)=1}}\chi_d(N)^{|J|}=\begin{cases}
    0,&\hbox{if $2\nmid|J|$,}\\
    1,&\hbox{if $2\mid|J|$.}    
\end{cases}$$
We thus only retain the terms with $2\mid |J|$, obtaining
$$
\begin{aligned}
\mathcal{O}_{E,k}(X;S;M)\approx\sum_{\substack{J\subseteq\{1,\dots,k\}\\2\mid\ |J|}}&\sumstar_{\substack{d\geq 1 \\ (d,N)=1}}g\bfrac dX\chi_d(M)\prod_{j\in J}\frac{\mathcal{X}_E(s_j)}{d^{2s_j-1}}\\
&\cdot\quad\sum_{n_1,\dots,n_k\geq 1}\frac{\lambda_E(n_1)\dots \lambda_E(n_k)\chi_d(n_1\dots n_k)}{n_1^{s_1}\dots n_k^{s_k}}.
\end{aligned}
$$

Now a procedure similar to our estimation of \eqref{eqn: sjj} leads to
\begin{conjecture}[CFKRS~\cite{cfkrs}]\label{Recipe conjecture elliptic} Let $\mathcal{O}_{E,k}(X;S;M)$ be as defined in \eqref{eqn: ellipticmoment}. If $\lab\re(s)-\frac{1}{2}\rab\ll 1/\log X$ for all $s\in S$, then as $X\rightarrow \infty$ we have
    \begin{equation}\label{eqn: elliptic recipe}
    \mathcal{O}_{E,k}(X;S;M)\sim\sum_{\substack{J\subseteq\{1,\dots,k\}\\2\mid\ |J|}}\prod_{j\in J} \mathcal{X}_E(s_j)\cdot X^{1+|J|-2\sum\limits_{j\in J}s_j}\tilde g\lz1-2\sum_{j\in J}s_j+|J|\pz H_{E,M}(S\setminus S_J\cup S_J^{-}),
    \end{equation} where $S_J=\{s_j:j\in J\}$, $S_J^{-}=\{1-s_j :j\in J\}$, and

    \begin{equation}\label{eqn: HEMSdefn}
        H_{E,M}(S)=\frac{1}{2\zeta(2)}\sum_{Mn_1\dots n_k=\square}\frac{\lambda_E(n_1)\dots \lambda_E(n_k)a(MNn_1\dots n_k)}{n_1^{s_1}\dots n_k^{s_k}}
    \end{equation}
    for $S=\{s_1,\dots,s_k\}$, with $a(n)$ defined by \eqref{eqn: a(n)def}.
\end{conjecture}

For $k=1$, the above conjecture asserts that
$$
\sumstar_{(d,N)=1}g\bfrac dX\chi_d(M)L(s,E_d)\sim X\tilde g(1) \frac{1}{2\zeta(2)}\sum_{Mn=\square}\frac{\lambda_E(n)a(MNn)}{n^{s}}.
$$
However, comparing this with the asymptotic formula \eqref{eqn: appendix asymptotic for elliptic moment} below, we see that this is false because it is missing one term.
 
We discarded some terms due to the averaging over the root numbers. However, we observe that the root numbers correlate with the Dirichlet coefficients of $L(s,E_d)$. With this observation, we now carry out the second version of the recipe, as described earlier in this section. In this second version, we average the product of the root numbers and the Dirichlet series coefficients of $\prod_{j=1}^k L(s,E_d)$ together. This modified recipe prediction is consistent with the multiple Dirichlet series prediction in Section~\ref{sec: ellipticMDS} below.

We write \eqref{eqn: elliptic moment after approx fe} as
\begin{align*}
    \mathcal{O}_{E,k}(X;S;M)\approx\sum_{J\subseteq\{1,\dots,k\}}& \sumstar_{\substack{d\geq 1\\(d,N)=1}}g\bfrac dX\chi_d(M)\prod_{j\in J}\frac{\mathcal{X}_E(s_j)}{d^{2s_j-1}}\\
    &\cdot\quad \sum_{n_1,\dots,n_k\geq 1}\frac{\lambda_E(n_1)\dots \lambda_E(n_k)\chi_d(N^{|J|}n_1\dots n_k)}{n_1^{s_1}\dots n_k^{s_k}}.
\end{align*}
Then, approximating the $d$-sum as in \eqref{eqn: dsumaverage}, we make the following prediction, which is similar to Conjecture~\ref{Recipe conjecture elliptic} but includes the terms with $2\nmid|J|$ in the asymptotic formula.
\begin{conjecture}[CFKRS~\cite{cfkrs}]\label{conj: modified recipe elliptic} Let $\mathcal{O}_{E,k}(X;S;M)$ be as defined in \eqref{eqn: ellipticmoment}. If $\lab\re(s)-\frac{1}{2}\rab\ll 1/\log X$ for all $s\in S$, then as $X\rightarrow \infty$ we have
    \begin{equation}\label{eqn: elliptic recipe MDS}
    \begin{aligned}
    \mathcal{O}_{E,k}(X;S;M)\sim \sum_{J\subseteq\{1,\dots,k\}}&\prod_{j\in J} \mathcal{X}_E(s_j)\cdot X^{1+|J|-2\sum_{j\in J}s_j}\tilde g\lz1+|J|-2\sum_{j\in J}s_j\pz\\
    &\cdot H_{E,MN^{|J|}}(S\setminus S_J\cup S_J^{-}),
    \end{aligned}
    \end{equation}
    where $S_J=\{s_j:j\in J\}$, $S_J^{-}=\{1-s_j :j\in J\}$, and $H_{E,M}(S)$ is defined by \eqref{eqn: HEMSdefn}.
\end{conjecture}

\begin{remark}
\begin{enumerate}      
\item Note that
$$
H_{E,MN^{|J|}}(S)=\begin{cases}
    H_{E,M}(S),&\hbox{if $2\mid |J|$},\\
    H_{E,MN}(S),&\hbox{if $2\nmid |J|$}.
\end{cases}
$$
Thus, for each $J$ with $2\mid |J|$, the term corresponding to $J$ in \eqref{eqn: elliptic recipe MDS} is the same as the corresponding term in \eqref{eqn: elliptic recipe}.
\item For $k=1$, Conjecture~\ref{conj: modified recipe elliptic} predicts
$$\sumstar_{(d,N)=1}g\bfrac dX\chi_d(M)L(s,E_d)\sim X\tilde g(1) H_{E,M}(s)+X^{2-2s}\mathcal X_E(s)\tilde g(2-2s)H_{E,MN}(-s),$$ which agrees with \eqref{eqn: appendix asymptotic for elliptic moment} below.

\item CFKRS \cite[Section 4.4]{cfkrs} restrict the family to arithmetic progressions $d\equiv a\mod N$, which makes the root numbers constant throughout the family. Summing their conjecture over $a\mod N$ leads to Conjecture~\ref{conj: modified recipe elliptic}.
\end{enumerate}
\end{remark}

\subsection{Quadratic twists of elliptic curves: predictions via multiple Dirichlet series}\label{sec: ellipticMDS}
An application of Mellin inversion gives
\begin{equation}\label{eqn: elliptic integral}
\mathcal{O}_{E,k}(X;S;M)=\frac{1}{2\pi i}\int_{(c)}A_{E,M}(s_1,\dots,s_k;w)X^w\tilde g(w)dw,\end{equation}
where
\begin{equation}\label{eqn: AEMdefn}
A_{E,M}(s_1,\dots,s_k;w)=\sumstar_{\substack{d\geq1 \\(d,N)=1}}\frac{\chi_d(M)L(s_1,E_d)\dots L(s_k,E_d)}{d^w}.
\end{equation}
\begin{prop}\label{prop: AEMregion}
    The multiple Dirichlet series $A_{E,M}(s_1,\dots,s_k;w)$ is absolutely convergent for $\re(s_j)>1$, $\re(w)>1$. It has a meromorphic extension to the region $\re(w)>1/2$, $\re(s_j)>5/4$ with only a simple pole at $w=1$ with residue
    $$\res{w=1}A_{E,M}(s_1,\dots,s_k;w)=H_{E,M}(S).$$
\end{prop}
\begin{proof}
    The proof is similar to the proof of Proposition \ref{prop: symplectic region}.
\end{proof}

The contribution of the pole at $w=1$ to \eqref{eqn: elliptic integral} is
$$X\tilde g(1)H_{E,M}(S),$$ which corresponds to the diagonal contribution $J=\emptyset$ in \eqref{eqn: elliptic recipe}.

The other terms come from the functional equations. Let $J\subseteq\{1,\dots,k\}$. Using the functional equation \eqref{eqn: elliptic twist functional equation} for $L(s_j,E_d)$ with $j\in J$, we deduce from \eqref{eqn: AEMdefn} that
\begin{align}
A_{E,M}(s_1,\dots,s_k;w)
&=\prod_{j\in J}\mathcal{X}_E(s_j)\sumstar_{\substack{d\geq1 \\(d,N)=1}}\frac{\chi_d(M)\chi_d(N)^{|J|}L(s_1^J,E_d)\dots L(s_k^J,E_d)}{d^{w+|J|-2\sum_{j\in J}s_j}} \notag\\
&=\prod_{j\in J}\mathcal{X}_E(s_j)\cdot A_{E,MN^{|J|}}\lz s_1^J,\dots,s_k^J;w+2\sum_{j\in J}s_j-|J|\pz. \label{eqn: AEMfunceqn}
\end{align}
Let $$\sigma_J(w):=w+2\sum_{j\in J}s_j-|J|.$$
The functional equation \eqref{eqn: AEMfunceqn} gives rise to a new pole at
$$
\sigma_J(w)=1,
$$
with residue
$$
\prod_{j\in J}\mathcal{X}_E(s_j)\cdot\res{w=1}A_{E,MN^{|J|}}(s_1^{J},\dots,s_k^J;w)=\prod_{j\in J}\mathcal{X}_E(s_j)\cdot H_{E,MN^{|J|}}(S\setminus S_J\cup S_J^-)
$$
by Proposition~\ref{prop: AEMregion}. The contribution of this residue to the integral \eqref{eqn: elliptic integral} is
$$
\prod_{j\in J}\mathcal{X}_E(s_j)\cdot X^{1+|J|-2\sum_{j\in J}s_j}\tilde g\lz 1+|J|-2\sum_{j\in J}s_j\pz H_{E,MN^{|J|}}(S\setminus S_J\cup S_J^-)
$$

As in the earlier sections, this conjecture is a consequence of the following natural conjecture about the meromorphic continuation of $A_{E,M}(s_1,\dots,s_k;w)$
\begin{conjecture}\label{conj: continuation elliptic}
    The function $$\prod_{J\subseteq\{1,\dots,k\}}(\sigma_J(w)-1)\cdot A_{E,M}(s_1,\dots,s_k;w)$$ has a holomorphic continuation to a tube domain that contains the point
    $$
    (s_1,\dots,s_k,w) = (\tfrac{1}{2},\dots,\tfrac{1}{2},1),
    $$
    and it is polynomially bounded in this region.
\end{conjecture}
Similarly as in the previous sections, we have the following Proposition.
\begin{prop}\label{prop: MDSelliptic}
    Conjecture \ref{conj: continuation elliptic} implies Conjecture \ref{conj: modified recipe elliptic} with a power-saving error term.
\end{prop}

\section{Appendix: the first moment of quadratic twists of elliptic curves}
Our goal here is to compute the asymptotic formula with power-saving error term for the first moment of the family of quadratic twists of the $L$-function of an arbitrary fixed elliptic curve over $\mathbb{Q}$. 

There is a number of papers in the literature that consider moments of the family of quadratic twists of $L$-functions associated to modular forms. Most of these papers focus on the case of forms of full level and note that the result can be extended to a general congruence subgroup. We explain here the necessary modifications in the case of the first moment. The results apply to quadratic twists of elliptic curve $L$-functions because of the modularity theorem due to Breuil, Conrad, Diamond, Taylor, and Wiles~\cite{bcdt,taylorwiles,wiles}.

We first give an overview of the proof of the asymptotic formula of the first moment in the case of forms of full level. This proof is due to Shen~\cite{shen}. Afterwards, we explain the needed modifications to deal with quadratic twists of $L(s,E)$.

\subsection{The level one case}
Shen~\cite{shen} considers the moment
$$M(\alpha,\ell):=\sum_{d\geq 1} \mu^2(2d) g\bfrac dX \chi_{8d}(\ell)L(\tfrac{1}{2}+\alpha,f\otimes\chi_{8d}),$$ where $f$ is a Hecke eigenform of weight $\kappa$ for the full modular group $SL_2(\Z)$ and $g$ is a smooth, compactly supported test function. He proves that for the shift $\alpha$ in a certain range,
\begin{equation}\label{eqn: main term of level 1}
M(\alpha,\ell)=X\tilde g (1) R(\alpha,\ell)+X^{1-2\alpha}  i^\kappa \tilde g(1-2\alpha)\bfrac{8}{2\pi}^{-2\alpha}\frac{\Gamma\lz 1-\alpha\pz}{\Gamma\lz 1+\alpha\pz} R(-\alpha,\ell)+O(\ell^{1/2+\epsilon}X^{1/2+\epsilon}),
\end{equation} where the main term $R(\alpha,\ell)$ is defined in \cite[Conjecture 1.3]{shen}, and the two terms on the right-hand side of \eqref{eqn: main term of level 1} correspond to the 0- and 1-swap terms.

The proof begins by employing the approximate functional equation
$$L(\tfrac{1}{2}+\alpha,f\otimes\chi_d)=\sum_{n\geq 1}\frac{\lambda_f(n)\chi_d(n)}{n^{1/2+\alpha}}\omega_\alpha\bfrac{n}{d}+i^\kappa \bfrac{d}{2\pi}^{-2\alpha}\frac{\Gamma\lz\frac\kappa 2-\alpha\pz}{\Gamma\lz\frac\kappa 2+\alpha\pz}\sum_{n\geq 1}\frac{\lambda_f(n)\chi_d(n)}{n^{1/2-\alpha}}\omega_{-\alpha}\bfrac{n}{d},$$
where $\omega_\alpha$, $\omega_{-\alpha}$ are smooth weights defined in \cite[Lemma 2.1]{shen}.
We can then write the moment as \cite[p.~321]{shen}
\begin{equation}\label{eqn: appendix moment decomposition}M(\alpha,\ell)=M^+(\alpha,\ell)+M^-(\alpha,\ell)
\end{equation}
with
\begin{equation}\label{eqn: appendix def of M^+ and M^-}
\begin{aligned}
    M^+(\alpha,\ell)&=\sum_{d\geq 1} \mu^2(2d) g\bfrac dX\sum_{n\geq 1}\frac{\lambda_f(n)\chi_{8d}(\ell n)}{n^{1/2+\alpha}}\omega_\alpha\bfrac{n}{8d}, \ \ \ \ \text{and}\\
    M^-(\alpha,\ell)&=i^\kappa \sum_{d\geq 1} \mu^2(2d) g\bfrac dX \bfrac{8d}{2\pi}^{-2\alpha}\frac{\Gamma\lz\frac\kappa 2-\alpha\pz}{\Gamma\lz\frac\kappa 2+\alpha\pz}\sum_{n\geq 1}\frac{\lambda_f(n)\chi_{8d}(\ell n)}{n^{1/2-\alpha}}\omega_{-\alpha}\bfrac{n}{8d}.
\end{aligned}
\end{equation}
We have
\begin{equation}\label{eqn: appendix M^- in terms of M^+}
\begin{aligned}
M^{-}(\alpha,\ell)&=X^{-2\alpha} i^\kappa \bfrac{8}{2\pi}^{-2\alpha} \frac{\Gamma\lz\frac\kappa 2-\alpha\pz}{\Gamma\lz\frac\kappa 2+\alpha\pz} \sum_{d\geq 1} \mu^2(2d) g\bfrac dX\bfrac{d}{X}^{-2\alpha} \\
&\cdot\quad\sum_{n\geq 1}\frac{\lambda_f(n)\chi_{8d}(\ell n)}{n^{1/2-\alpha}}\omega_{-\alpha}\bfrac{n}{8d},
\end{aligned}
\end{equation}
and we can note that the $d,n$-sum is $M^+(-\alpha,\ell)$ with $g(x)x^{-2\alpha}$ in place of $g(x)$.
Shen then proceeds to prove that (see equations (3.10) and (7.3) and Lemmas~5.3 and 6.1 of \cite{shen})
\begin{equation}\label{eqn: appendix asymptotic for M^+}
M^+(\alpha,\ell) =  X\tilde g(1) R(\alpha,\ell) +O(\ell^{\varepsilon} X^{\frac{1}{2}+\varepsilon})
\end{equation} and the result follows after combining \eqref{eqn: appendix moment decomposition}, \eqref{eqn: appendix def of M^+ and M^-}, \eqref{eqn: appendix M^- in terms of M^+} and \eqref{eqn: appendix asymptotic for M^+}.

\subsection{The elliptic curves case}
Now we consider an elliptic curve $E$ over $\mathbb{Q}$ of conductor $N$, and the family of quadratic twists $E_d$, where $d$ runs over positive fundamental discriminants coprime to $N$. 

We are interested in moments of the form
$$M_E(\alpha,\ell)=\sum_{\substack{d\geq 1\\(d,N)=1}} \mu^2(2d)g\bfrac{d}{X}\chi_{8d}(\ell) L(\tfrac{1}{2}+\alpha,E_{8d}),$$
where the coprimality condition can also be captured by ensuring that $N^2|\ell$.

The main difference from the level one case comes from the approximate functional equation.
Accounting for the sign in the functional equation \eqref{eqn: elliptic twist functional equation}, we have (recall that now the weight $\kappa =2$)
$$
\begin{aligned}L(\tfrac{1}{2}+\alpha,E_d)=&\sum_{n\geq 1}\frac{\lambda_E(n)\chi_d(n)}{n^{1/2+\alpha}}\omega_\alpha\bfrac{n}{d}\\
& + \varepsilon(E)\chi_d(N) \bfrac{d  \sqrt N }{2\pi}^{-2\alpha}\frac{\Gamma\lz 1 -\alpha\pz}{\Gamma\lz 1 +\alpha\pz}\sum_{n\geq 1}\frac{\lambda_E(n)\chi_d(n)}{n^{1/2-\alpha}}\omega_{-\alpha}\bfrac{n}{d},
\end{aligned}$$ which leads to
$$M_E(\alpha,\ell)=M_E^+(\alpha,\ell)+M_E^-(\alpha,\ell)$$ with
$$\begin{aligned}
    M_E^+(\alpha,\ell)&=\sum_{\substack{d\geq 1\\(d,N)=1}} \mu^2(2d) g\bfrac dX\sum_{n\geq 1}\frac{\lambda_E(n)\chi_{8d}(\ell n)}{n^{1/2+\alpha}}\omega_\alpha\bfrac{n}{8d}, \ \ \ \ \text{and}\\
    M_E^-(\alpha,\ell)&= N^{-\alpha} \varepsilon(E)\sum_{\substack{d\geq 1\\(d,N)=1}} \mu^2(2d) g\bfrac dX \bfrac{8d}{2\pi}^{-2\alpha}\frac{\Gamma\lz 1-\alpha\pz}{\Gamma\lz 1 +\alpha\pz}\sum_{n\geq 1}\frac{\lambda_E(n)\chi_{8d}(\ell nN)}{n^{1/2-\alpha}}\omega_{-\alpha}\bfrac{n}{8d}.
\end{aligned}$$
We now see that $$
\begin{aligned}
M_E^-(\alpha,\ell)=& \varepsilon(E)X^{-2\alpha} \bfrac{8  \sqrt N }{2\pi}^{-2\alpha}\frac{\Gamma\lz 1-\alpha\pz}{\Gamma\lz 1+\alpha\pz}\sum_{\substack{d\geq 1\\(d,N)=1}} \mu^2(2d) g\bfrac dX \bfrac{d}{X}^{-2\alpha}\\
&\cdot\quad\sum_{n\geq 1}\frac{\lambda_E(n)\chi_{8d}(\ell nN)}{n^{1/2-\alpha}}\omega_{-\alpha}\bfrac{n}{8d},
\end{aligned}$$ where the $d,n$-sum is now $M_E^+(-\alpha,\ell N)$ with $g(x)x^{-2\alpha}$ in place of $g(x)$. A straightforward modification of the proof of Shen now gives an asymptotic formula
\begin{equation}\label{eqn: appendix asymptotic for elliptic moment}M_E(\alpha,\ell)\sim X\tilde g(1)R_E(\alpha,\ell)  + \varepsilon(E)X^{1-2\alpha} \tilde g(1-2\alpha)   \bfrac{8\sqrt N}{2\pi}^{-2\alpha}  \frac{\Gamma\lz 1-\alpha\pz}{\Gamma\lz 1+\alpha\pz} R_E(  - \alpha,\ell N)\end{equation} with a power-saving error term. We thus see that the result again contains both the 0- and 1-swap terms.

\end{document}